\documentclass[manuscript,screen]{acmart}

\setcopyright{rightsretained}
\usepackage{xcolor}
\usepackage{maple}
\definecolor{mygreen}{rgb}{0,0.6,0}
\definecolor{mygray}{rgb}{0.5,0.5,0.5}
\definecolor{mymauve}{rgb}{0.58,0,0.82}
\definecolor{altblue}{rgb}{0.0,0.6,1.0}
\definecolor{lstbg}{cmyk}{0.05, 0.01, 0, 0}
\definecolor{morebluish}{cmyk}{0.06,0.04,0,0}
\usepackage{listings}
\input{listings-maple-definition.sty}
\newcommand{\Ce}{\ensuremath{\mathrm{Ce}}}
\newcommand{\ce}{\ensuremath{\mathrm{ce}}}

\newcommand{\e}{\ensuremath{\varepsilon}}
\newtheorem{remark}{Remark}[section]
\begin{document}

\title{Blendstrings:\\ an environment for computing with smooth functions}

\author{Robert M. Corless}
\email{rcorless@uwo.ca}
\orcid{0000-0003-0515-1572}
\affiliation{%
  \institution{Editor-in-Chief, Maple Transactions}
  \city{London}
  \state{Ontario}
  \country{Canada}
  \postcode{N6A 5B9}
}


\begin{abstract}
A ``blendstring'' is a piecewise polynomial interpolant with high-degree two-point Hermite interpolational polynomials on each piece, analogous to a cubic spline. Blendstrings are smoother and can be more accurate than cubic splines, and can be used to represent smooth functions on a line segment or polygonal path in the complex plane. I sketch some properties of blendstrings, including efficient methods for evaluation, differentiation, and integration, as well as a prototype Maple implementation. Blendstrings can be differentiated and integrated exactly and can be combined algebraically. I also show applications of blendstrings to solving differential equations and computing Mathieu functions and generalized Mathieu eigenfunctions.
\end{abstract}



\received{31 January 2023}

\maketitle

\section{Introduction}
Providing useful tools for dealing with mathematical functions is one of the fundamental tasks of symbolic computation systems. This usually entails various kinds of polynomial or rational approximations, which are used to provide accurate approximations of the values of the function and its derivatives and antiderivatives.  The various major systems use different methods, and may even use different methods (hidden from the user) for different functions in the same system.  Arbitrary-precision evaluation of the exponential function, for instance, might use argument reduction, rational approximation, and powering. A common model is that the user treats the system's representation as an oracle, and asks for evaluation of the function at various points.  The system may or may not use any previous evaluation points as a shortcut in evaluating at the new point.

In contrast, the Chebfun system~\cite{Battles2004} uses piecewise Chebyshev polynomials with fast construction methods based on the FFT, and is very efficient and very reliable.  The basic model in the Chebfun system, however, is more like simultaneous evaluation at a great many points, although rapid change to Chebyshev series via the FFT is also fundamental.  Combination of these functions then becomes a question of numerical combinations of the points.  This can be orders of magnitude faster than traditional computer algebra systems can combine or manipulate their typical representations.  The ApproxFun system makes even further improvements in speed, accuracy, and generality; see~\url{https://juliaapproximation.github.io/ApproxFun.jl/latest/}.

This paper describes something similar to Chebfun and ApproxFun, but which uses what I call ``blendstrings.'' A formal definition will be given in section~\ref{sec:Definitions}, but, in brief, blendstrings are piecewise polynomials where each polynomial is represented not in any standard basis but rather as a \textsl{blend}. A blend is a two-point Hermite interpolant, typically of quite high degree~\cite{Corless2021}.  To construct a blendstring, one needs Taylor coefficients at several knots, typically lying in a straight line in the complex plane and even more typically on an interval of the  real line.

The basic method used to evaluate a blend is described in~\cite{Corless2021}.
The reasons that blends are potentially interesting objects to work with are laid out in~\cite{Corless2023details}.  This present paper shows how to use these objects in a composite fashion, \textsl{i.e.} as blendstrings, and demonstrates an application where this is useful. Using blends in this composite fashion is similar to the extension from a single cubic polynomial to a cubic spline.

Blendstrings are \textsl{not} as good as chebfuns or approxfuns in many ways.  They do not have quite as good approximation properties, for instance, which can be seen by their Lebesgue constants, which are $O(\sqrt{n})$ in size whereas those of Chebyshev series are $O(\ln n)$, as the degree $n\to \infty$.  Nonetheless blends seem natural in some contexts, such as the numerical solution of Ordinary Differential Equations (ODEs) by high-order marching methods~\cite{nedialkov2005solving}. Such marching methods are of cost polynomial in the number of bits of accuracy requested~\cite{ilie2008adaptivity}, and are especially interesting for functions whose Taylor coefficients can be computed easily at a given point, i.e. D-finite or holonomic functions~\cite{vanderHoeven1999,van2001fast,mezzarobba2012note}. Blendstrings have an arbitrarily high order of continuity, depending only on how many Taylor coefficients are known at each knot. They might be especially useful in arbitrary-precision environments.  They have been found to be useful in the study of Mathieu functions~\cite{Brimacombe2021}.

This paper describes a prototype implementation in Maple and some of its capabilities.  This implementation is available at~\url{https://github.com/rcorless/Blends-in-Maple/blob/main/BlendstringExamples.maple}.

The implementation uses \lstinline{evalhf} for speed when possible, but is primarily intended for high precision.

\section{Definitions and basic properties\label{sec:Definitions}}
In what follows, the word \textsl{grade} means ``degree at most''. That is, a polynomial of grade (say)~$5$ is of degree at most~$5$. But with blends one does not immediately know the exact degree because the leading coefficients are not visible and might in actuality be zero.

Consider an analytic function $f(z)$ with Taylor series coefficients known at $z=a$ and at $z=b$.
This paper assumes that the Taylor coefficients are known to a fixed working precision, most commonly $15$ decimal digits. If one wants to work in higher precision, it will be necessary to start with Taylor coefficients evaluated to that higher precision.
Convert to the unit interval by introducing a new variable~$s$ with $z = a + s(b-a)$.
\subsection{The basic formula}
The following formula, known already to Hermite~\cite[p.~4]{hermite1873cours}, has the property that the grade $m+n+1$ polynomial
\begin{align}
    H_{m,n}(s) = &\sum _{j=0}^{m} \left[ \sum _{k=0}^{m-j}{n+k\choose k}{s}^{k+j}
 \left( 1-s \right) ^{n+1}\right]p_{{j}}  \nonumber\\
 +& \sum _{j=0}^{n}
\left[\sum _{k=0}^{n-j}{m+k\choose k}{s}^{m+1}
 \left( 1-s \right) ^{k+j}\right] \left( -1 \right) ^{j}q_{{j}}  \label{eq:TPHI}
\end{align}
has a Taylor series matching the given $m+1$ values $p_j = f^{(j)}(0)/j!$ at $s=0$ and another Taylor series matching the given $n+1$ values $q_j = f^{(j)}(1)/j!$ at $s=1$.  Putting this in symbolic terms and using a superscript $(j)$ to mean the $j$th derivative with respect to~$s$ gives
\[
\frac{H^{(j)}_{m,n}(0)}{j!} = p_j \qquad \mathrm{and} \qquad \frac{H^{(j)}_{m,n}(1)}{j!} = q_j
\]
for $ 0 \le j \le m$ on the left and for
$0 \le j \le n$ on the right.
This is a kind of interpolation, indeed a special case of what is called \emph{Hermite} interpolation.
As with Lagrange interpolation, where for instance two points give a grade one polynomial, that is, a line, here $m+n+2$ pieces of information gives a grade $m+n+1$ polynomial.  As per~\cite{Corless2021} this formula can be evaluated in $O(m+n)$ arithmetic operations.

\subsection{Definition of a blendstring}
A \textsl{blendstring} is an ordered finite set of ``local Taylor polynomials''
\[
\mathrm{LTP}_k(z) := c_{k,0} + c_{k,1}(z-a_k) + c_{k,2}(z-a_k)^2 + \cdots + c_{k,m_k}(z-a_k)^{m_k}
\]
for $0 \le k \le M$,
together with the line segments $[a_0,a_1]$, $[a_1, a_2]$, $\ldots$, $[a_{m_k-1}, a_{m_k}]$, over which one can blend the appropriate Taylor polynomials.  The grades $m_k\ge 0$ of the Taylor polynomials are integers.  It is not necessary that the leading coefficients $c_{k,m_k}$ be nonzero. It is necessary that the knots $a_k$ be distinct from their predecessor knot and their successor knot (but crossings are permitted otherwise: polygonal paths in the complex plane can return to earlier knots). The local Taylor polynomials may be represented by an array with the knot $a_k$ first and then all the Taylor coefficients: $[a_k, c_{k,0}, \ldots, c_{k,m_k}]$.

Two blendstrings
\[
\mathcal{B}_1 = \{[a_k, c_{k,0}, \ldots, c_{k,m_k}]\}_{k=0}^M
\]
and
\[
\mathcal{B}_2 = \{[b_k, d_{k,0}, \ldots, d_{k,n_k}]\}_{k=0}^N
\]
are \textsl{compatible} if $N=M$ and all corresponding knots $a_k = b_k$ and all corresponding grades $m_k = n_k$.
\subsection{Approximation theoretic properties}
Recall that the Lebesgue function $L(s)$ for a polynomial basis $\phi_j(s)$ for $0 \le j \le N$ is the function $L(s) = \sum_{j=0}^N |\phi_j(s)|$.  Relative errors $\delta_j$ in polynomial coefficients $a_j(1+\delta_j)$ produce changes $\Delta p(s)$ in the polynomial value. By the triangle inequality, $|\Delta p(s)| \le L(s) \|\mathbf{a}\|_\infty \varepsilon$ if all relative coefficient changes $|\delta_j| \le \varepsilon$.
It was claimed in~\cite{Corless2021} that for a \textsl{balanced} blend, that is one where the grade $m$ of Taylor polynomial at the left end is the same as the grade $n$ at the right end, $L(s) \le 2$ on $0 \le s \le 1$ so that the Lebesgue constant, on this interval, is just $2$, independently of the degree. On the more usual interval of $[-1,1]$, the Lebesgue constant $\Lambda_{m,m} \sim 2\sqrt{m/\pi}$ which is not much worse than the optimal growth, which is $O(\ln m)$. The proofs of these claims can be found in~\cite{Corless2023details}.

We also claimed in~\cite{Corless2021} that a double Horner evaluation using IEEE 854 floating-point arithmetic was backward stable on $0 \le s \le 1$ in the following sense: namely, that numerical evaluation with unit roundoff $u$ gave the \textsl{exact} answer for a blend with Taylor coefficients $p_j$ and $q_j$  changed to $p_j(1+\delta_{p,j})$ and $q_j(1+\delta_{q,j})$ with $|\delta_{\cdot,j}| \le \gamma_{3(m+n)}$ where $\gamma_k = ku/(1-ku)$. A proof can be found in~\cite{Corless2023details}.  This theorem, together with the bound on $L(s)$, guarantees accurate approximation.

Blends and blendstrings make sense even if the coefficients are exactly known in terms of symbols such as $\gamma$, $\pi$ and the like, but they lose much of their advantage in rapid and stable computation.  Assume henceforth that the Taylor series at each knot are approximated to working precision in floating-point complex arithmetic.  IEEE-754 double precision is frequently used, but Maple's arbitrary precision arithmetic is sometimes convenient as well.  Working at a fixed precision, evaluation of a blend is of cost $O(m+n)$ multiplications at that precision.

The truncation error in approximating $f(z)$ by the blend is given by (in the scaled variable $s$)
\begin{equation}
f(s) - H_{m,n}(s) = \frac{f^{(m+n+2)}(\theta)}{(m+n+2)!} s^{m+1}(s-1)^{n+1} \>. \label{eq:errorblend}
\end{equation}
Here $\theta \in (0,1)$ is otherwise unknown.  The maximum value of the polynomial $s^{m+1}(1-s)^{n+1}$ on the interval $(0,1)$ is attained at the point $s=(m+1)/(m+n+1)$ and is $(m+1)^{m+1}(n+1)^{n+1}/(m+n+2)^{m+n+2}$, which if $m=n$ reduces to $2^{-2(m+1)}$.  One frequently uses $m=10$ or more, and so this factor reduces the error by a factor of a million or more, compared to a Taylor series expansion of grade $m+n+2$ on only one side.  In some sense this factor is the real reason blends are interesting.

Approximation by a blendstring is analogous to approximation by cubic splines, though of higher order.  If each subinterval of the blendstring is of width $h$, and each Taylor polynomial is of grade $m$, then the order of accuracy of the blendstring is $\mathcal{B}(z)-f(z) = O(h^{2m+2})$. More, the accuracy of the derivative is $\mathcal{B}'(z)-f'(z) = O(h^{2m+1})$, and so on, losing one order of accuracy per derivative.  The proof follows by standard methods and can be seen using equation~\eqref{eq:errorblend} once the polynomial $(z-a_n)^{m+1}(z-a_{n+1})^{m+1}$ is expressed using the transformation $z=a_n + sh$.  The quite high-order Taylor coefficients that appear will not cause problems for the smooth functions that are considered here, but can prove troublesome if there are nearby singularities, as expected\footnote{The detection and location of nearby singularities is an interesting problem in practice and an active topic of research. The best-known ways use information gleaned from local Taylor series, for instance by converting to Pad\'e approximants.}.
\subsection{Evaluation and differentiation of a blend}
The paper~\cite{Corless2021} implemented a double Horner expansion of the formula~\eqref{eq:TPHI}.   Several examples were given showing its numerical stability in practice. As stated earlier, the cost of evaluation is linear in the grade of the blend, $O(m+n)$.
Since the Horner-like evaluation is simply a pair of for-loops, automatic differentiation is straightforward. By experiment, rounding errors in the automatic derivatives are similarly small, even though differentiation is infinitely ill-conditioned.

\section{The prototype implementation}
The code can be found in the files \lstinline{Blend.mpl}, \lstinline{deval.mpl}, and \\
\lstinline{BlendstringUtilities.mpl}, at \url{https://github.com/rcorless/Blends-in-Maple/blob/main/BlendstringExamples.maple}.
\subsection{Constructing a blendstring}
The basic data structure I chose for the Maple implementation is a two-dimensional \lstinline{Array}.  Specifically, the following command gives an empty blendstring with space for $M+1$ knots and local Taylor polynomials of grade $m$:
\begin{lstlisting}
B := Array(0..M, 0..m+1 );
\end{lstlisting}
This limits each $m_k$ to be the same number $m$, but balanced blends\footnote{A \textsl{balanced blend} is one in which $m \approx n$. Exact equality is best, but near equality is almost as good. Highly unbalanced blends are exponentially bad~\cite{Corless2023details}.} are best anyway for approximation, and this sufficed for my first applications.  Provision for variable grades $m_k$ in the same blendstring should be done in a future version.
This structure represents the local Taylor polynomials densely in the style described earlier: knot first, coefficients later.  The most basic way of constructing a blendstring is simply to create such an object directly.

\begin{figure}
    \centering
    \includegraphics[width=0.25\textwidth]{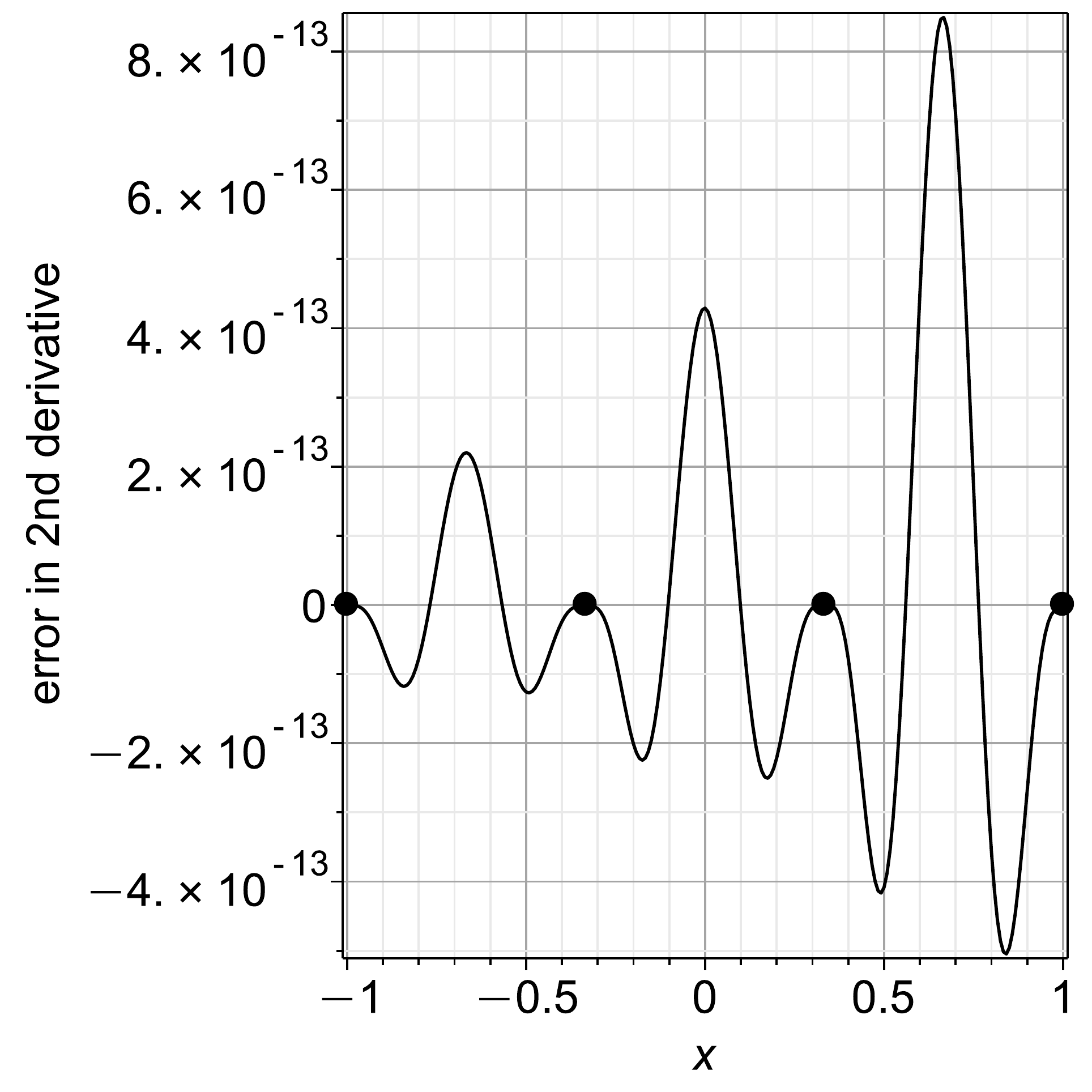}
    \caption{The error in the automatically computed second derivative of the blend for $\exp(z)$ just constructed with knots at $-1$, $-1/3$, $1/3$, and $1$, each with Taylor polynomials of grade $m=5$.}
    \label{fig:errblendexp2}
\end{figure}

Constructing a blendstring for a given analytic function $f(z)$ is then a matter of filling the \lstinline{Array} with the appropriate knots and Taylor coefficients, perhaps generated by the \lstinline{series} command. For example, the following constructs a blendstring with four knots and local Taylor polynomials of grade $5$.  This blendstring approximates $\exp(z)$ to better than $5\cdot 10^{-15}$ over the interval $(-1,1)$.  Figure~\ref{fig:errblendexp2} shows that the error in the 2nd derivative is smaller than $10^{-12}$.
\begin{lstlisting}
M := 3;
grade := 5;
Digits := 15;
knots := Array(0..M, k -> ( -1 + 2*k/M) ) ;
f := z->exp(z);
B := Array(0..M, 0..grade+1 );
for k from 0 to M do
  B[k,0] := knots[k];
  S := series( f(z), z=knots[k], grade+1 );
  for j from 0 to grade do
    B[k,j+1] := evalf(coeff(S,z-knots[k],j));
  end do;
end do:
\end{lstlisting}
\begin{remark}
The Maple environment variable \lstinline{Order}, or equivalently the third parameter in the call to \lstinline{series}, only sets the ``working grade'' for \lstinline{series} and does not guarantee that the results are actually that grade of Taylor polynomial.  An example is \lstinline{series( f(x), x, 6 )} with $f(x)=\sin(x)/x$.  The returned answer is $1 + O(x^4)$.  This difficulty will be ignored in this paper, apart from this mention.
\end{remark}
The next block uses commands that will be described in section~\ref{sec:deval}, to plot the error in the $2$nd derivative compared to the evaluation of the symbolic derivative.
\begin{lstlisting}
Digits := 30;
y := deval(B, nder = 3, nRefine = 80):
N := upperbound(y)[1]:
diffs := Array(0 .. N, 0 .. 4):
for k from 0 to N do
  diffs[k, 0] := y[k, 0];
  for j to 4 do
    diffs[k,j] := y[k,j]-(D@@(j-1))(f)(diffs[k,0]);
  end do;
end do:
X := Vector(N + 1, j -> y[j-1, 0]):
Y := Vector(N + 1, j -> diffs[j-1, 3]):
eplot := plot(X, Y, colour=black, symbol=point,
                     axes=boxed, gridlines=true);
\end{lstlisting}
\subsection{Arithmetic combinations of compatible blendstrings}
Since Taylor polynomials can be added together, it is a simple matter to form linear combinations $\alpha \mathcal{B}_1 + \beta\mathcal{B}_2$.  The result is a blend of the same grade as the components (even if some of the leading coefficients of the resulting Taylor polynomials are zero, that is still information about the function that gets used in the blend).

Multiplication requires the Cauchy product, so each Taylor coefficient of the product is
\[
p_{k,j} = \sum_{\ell=0}^j a_{k,\ell}b_{k,j-\ell}\>,
\]
for $0 \le j \le m_k$.
Unlike in the case of addition, however, this truncated product of two blends is not the blend of the exact product of the two underlying polynomials, because degrees increase with products; nonetheless it is an appropriate approximation to the product of the functions underlying the two blends being multiplied.

For the moment, \textsl{division} requires that each constant coefficient of the dividend be nonzero, because I do not yet have a method of blending Laurent series.  However, the underlying function of the dividend might have a zero on the interval, and this might lead to unexpected results.

These arithmetic combinations have been implemented in this package by the command \lstinline{zipBlendstrings}, which takes as input an arbitrary binary function $(a,b) \to f(a,b)$ and uses Maple's built-in \lstinline{series} command to carry out the necessary operations.  This is efficient (if the efficient routines in the \lstinline{PolynomialTools} package are used to convert back and forth from Maple polynomials) because \lstinline{series} is in the kernel.
Since \lstinline{series} is smart about indeterminate forms, this can be effective.

Simply by using this utility, one can build up a collection of useful blendstrings on a given sequence of knots, starting from the function $z$ which has the value $a_k$ at every knot $a_k$ and derivative $1$ and all higher derivatives zero.  One could then represent powers of $z$ and polynomials in $z$ as compatible blends by simply applying the operations in sequence. As a simple example, I used the Chebyshev recurrence relation to construct the Chebyshev polynomial $T_6(z)$ as a blendstring on the same knots as the above example.  It worked, and was perfectly accurate, which was unsurprising.

\begin{figure}
    \centering
    \includegraphics[width=0.305\textwidth]{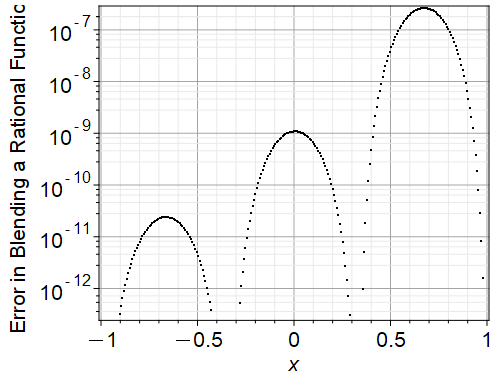}
    \caption{The absolute error $|\mathcal{B}(x_j) - (1+x_j/2)/(1-x_j/2)|$, computed at many points.  The blendstring $\mathcal{B}$ was constructed by constructing a blendstring for $1$ at four equally-spaced knots on the interval $[-1,1]$, a compatible one for $z$, a third for $1+z/2$, a fourth for $1-z/2$, and finally by dividing the one blendstring by the other.}
    \label{fig:pad11error}
    \Description{curves made up of black dots, on a gray grid background, making three rising equal-width ``humps" from left to right.}
\end{figure}

A more interesting test occurs when one computes a blendstring on those same knots, $\{-1, -1/3, 1/3, 1\}$, with all Taylor polynomials of grade $m=5$, for the \textsl{rational} function $(1+z/2)/(1-z/2)$.  This example tries to approximate a rational function on this interval with a piecewise polynomial, where the grade of the polynomial on each of the three subintervals is $m+n+1 = 11$.  The error in this approximation, which is appreciable because rational functions are not well approximated by polynomials, is shown in figure~\ref{fig:pad11error}.
\subsection{Applying a function to a blendstring to get another blendstring}
The package contains another function, \lstinline{mapOntoBlendstring}, which computes $f(\mathcal{B})$ given an operator representation\footnote{A Maple operator representation for a function $f(z)$ is \lstinline{z -> f(z)}.} for the function $f(z)$.  Again, this uses \lstinline{series} to compute local Taylor polynomials for $f( LTP_k(z) )$ about the knots $a_k$, and gives an approximation to the composition of functions.
\subsection{Evaluation and differentiation \label{sec:deval}}
Given the capability for efficient and numerically stable evaluation of a blend, all that remains is to decide how to organize the dispatch.  I chose not to use something like \lstinline{piecewise}, which is the standard way in Maple to represent a piecewise polynomial; instead I chose to model my evaluation of blendstrings on the subroutine \lstinline{deval} of the Matlab ODE Suite~\cite{shampine1997matlab}. The main reason for this was speed: by thinking of a blendstring as an ``all at once" object, with many function and derivative values already present, one can gain significant speedup.  The second reason was that I had erroneously thought that \lstinline{piecewise} was limited to the real line, and the application driving this implementation needs paths in the complex plane. However, that was a failure of imagination on my part, because \lstinline{piecewise} can indeed be used with complex input, so long as the boolean decision functions are correct.  But even so, the dispatch model seems geared to ``one at a time" inputs, whereas I really did want to ``think all at once," which meant to me something like \lstinline{deval}.  The \lstinline{deval} routine itself in Matlab is reminiscent of \lstinline{pp} (for ``piecewise polynomial") constructs in Matlab, which themselves tend to work ``all at once."

The call is either \lstinline{deval(y,'x'=pt)} to evaluate the blendstring \lstinline{y} at a single point \lstinline{pt}, or simply \lstinline{deval(y)} to evaluate the blendstring at (by default \lstinline{nRefine=2*Order}) equally-spaced points on each subinterval.  It is frequently the case that one wants to graph an entire blendstring, and using that many interior points is generally sufficient.  If not, one can ask for more (or fewer) by the \lstinline{nRefine=n} option, where $n$ is however many you want.

One reason for evaluating at so many points is for plotting. The current method for plotting is simply to hand the evaluated points to Maple's \lstinline{plot}; this can be inefficient if done with too many conversions back and forth from lists, so it would be well to write a specialized plot routine for the package.  This has not yet been done.

To differentiate, one uses the keyword \lstinline{nder}, as in \lstinline{deval(y,'nder'=2)} to compute (at all default interior points, as above) $y(x_k)$, $y'(x_k)$, and $y''(x_k)$.  The result is a two-dimensional \lstinline{Array}, with knots at all the original knots plus the new evaluations, together with derivatives (not Taylor coefficients) to all requested orders.

\textbf{An important caveat}.  Blends provide good approximations only on the segment $[a_k,a_{k+1}]$ which gets mapped to $0 \le s \le 1$, or very near in the complex plane to that segment.  Away from that segment, the truncation error grows very rapidly, like $|s|^{m+n+1}$.  So if one is using a blendstring to approximate a function in the complex plane, but not actually \textsl{on} the segment, getting accurate answers is not easy and may not be possible.  Further, if the segments of the blendstring follow a polygonal path or even make loops, it may not be evident which piece of the blendstring would be best to use to evaluate the function with.  In practice I have only used the \lstinline{nRefine} option in cases such as this, which guarantees that each evaluation point is precisely located exactly on one of the defining segments.
\subsection{Integration}
It is integration for which blends and blendstrings truly shine.  There is a formula reported in~\cite{Corless2021} for integrating a blend.
\begin{align}
\int_{s=0}^1 H_{m,n}(s)\,ds = &
{\frac { \left( m+1 \right) !}{ \left( m+n+2 \right) !}\sum _{j=0}^{m}
{\frac { \left( n+m-j+1 \right) !}{ \left( j+1 \right)
 \left( m-j \right) !}}}\,p_{{j}}\nonumber\\
 &+{\frac { \left( n+1 \right) !}{ \left( m+n+2
 \right) !}\sum _{j=0}^{n}{\frac {  \left( n+m-j
+1 \right) !}{ \left( j+1 \right)  \left( n-j \right) !}}}\,\left( -1 \right) ^{j}q_{{j}}
\>.
\label{eq:sumasintegral}
\end{align}
This is an exact complete integral across the subinterval, if the coefficients are known exactly, but the main use of this routine is when the coefficients are floating-point numbers, in which this becomes a kind of numerical quadrature. From this formula, one can construct an ``{exact}'' blend for the indefinite integral across the blend because now the Taylor coefficients \textsl{for the integral} are known at each end.  It is a simple matter to propagate constants along from one end of a blendstring to the other to construct an ``exact'' indefinite integral of the original blendstring.


This is an exact integral \textsl{for the blendstring}.  This gives an \textsl{approximation} to the integral of the underlying function $f(z)$ that the blendstring approximates.  Theory predicts that the results will be most satisfactory for balanced blends~\cite{Corless2023details}.

\subsubsection*{Accuracy and stability}
Because the formula is principally of use with floating-point approximations to the Taylor coefficients, the following theorem is useful.  Let the bold symbol $\mathbf{1}$ represent the series with all coefficients equal to $1$, and the bold symbol $\mathbf{(-1)^k}$ represent the series with all coefficients alternating in sign starting with $(-1)^0 = 1$.

\begin{theorem}
If the coefficients of the blend are in error by at most $\Delta p_j$ for $0 \le j \le m$ and $\Delta q_j$ for $0 \le j \le n$, then the error in the integral of the blend is bounded by $\int_0^1 H_{m,n}(\mathbf{1},\mathbf{(-1)^k}) \max |\Delta p_j|, |\Delta q_j|$.
\end{theorem}
The proof is immediate by using the linearity of the blend and the linearity of the integral.  Further, the integral can be explicitly computed:
\begin{equation}
\int_0^1 H_{m,n}(\mathbf{1},\mathbf{(-1)^k}) = 2\Psi(n+m+3) - \Psi(m+3) - \Psi(n+3) + \frac{n+m+4}{(n+2)(m+2)}\>.
\end{equation}
Here $\Psi(n+1) = -\gamma + \sum_{k=1}^n 1/k$ is the logarithmic derivative of the factorial function.  If either $n$ or $m$ goes to infinity while the other remains fixed, this integral grows like $\ln n$ or $\ln m$.  If both $n=m$ go infinity together the integral is asymptotic to $2\ln2 - 1/(2m) + O(1/m^2)$.

This shows that for balanced blends, the computation of the integral by using this formula is numerically stable, and that is certainly what is observed in practice.

\subsubsection*{An example}
\begin{figure}
    \centering
    \includegraphics[width=0.25\textwidth]{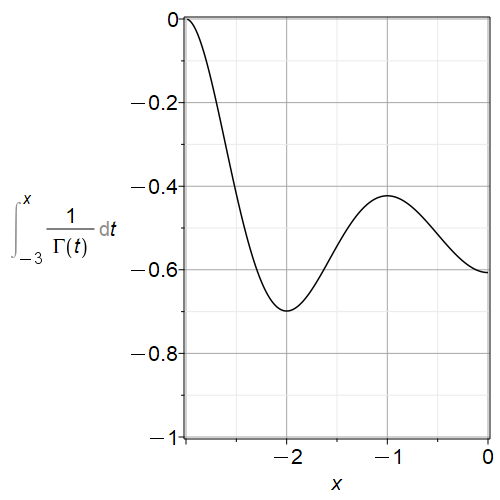}
    \caption{A plot of a blendstring for the indefinite integral of $1/\Gamma(x)$.}
    \label{fig:intrecipgamma}
\end{figure}

Here is an example, a quadrature for the function $1/\Gamma(s)$ on the interval $-3 \le s \le 0$.  The code below uses a blend of the known series at the negative integers, and grade $7$ Taylor polynomials at each knot:
\begin{lstlisting}
Digits := 15;
knots := [seq(-3 + k, k = 0 .. 3)]:
grade := 7:
B := Array(0 .. numelems(knots) - 1, 0 .. grade+1):
for k from 0 to numelems(knots) - 1 do
  B[k, 0] := knots[k + 1];
  S := series(1/GAMMA(s),s=knots[k + 1],grade+1);
  for j to grade+1 do
    B[k,j] := evalf(coeftayl(S,s=knots[k+1],j-1));
  end do;
end do:
Y := CodeTools:-Usage( intBlend(B) ):
CodeTools:-Usage(evalf(Int(1/GAMMA(s),s=-3..0)));
evalf(Y[numelems(knots) - 1, 1]);
\end{lstlisting}
Maple's built-in numerical integrator~\cite{geddes1992hybrid} gives $-0.606607588776539$ while the blend integration gives $-0.6066075887\textcolor{red}{83124}$.  Using grade $8$ Taylor polynomials instead gives better agreement, and using grade $10$ agrees to $17$ decimals if one works in $30$ Digits.  It is not fair to compare the timing of these (\lstinline{intBlend} is much faster on this example), because the series coefficients of $1/\Gamma(s)$ are known symbolically at negative integers and essentially free to evaluate numerically, while the general-purpose built-in numerical integrator does not take advantage of that information for this specific problem.  On other examples the built-in integrator can be faster.  On the other hand, the integral blend above can be evaluated efficiently at many points between $-3$ and $0$; it has effectively computed the \textsl{indefinite} integral.  This could be advantageous for some situations.
\begin{lstlisting}
y := CodeTools:-Usage(deval(Y, nRefine = 80)):
npts := upperbound(y)[1]:
\end{lstlisting}
The evaluation of the blendstring for the antiderivative of $1/\Gamma(x)$ took 31ms and thereafter the plot, shown in Figure~\ref{fig:intrecipgamma}, took no measureable time.  In contrast, the plot of the built-in version took nearly two seconds on the same machine.
\section{Applications}
There are many methods for approximating functions, and many applications for approximation of functions~\cite{trefethen2019approximation}.  The real question is, are there any applications for which blendstrings might be especially suited?  Their main benefits are a high degree of smoothness: if all the grades of Taylor coefficients are $m$, then the resulting blendstring is in $\mathcal{C}^m$. They require, however, derivatives (Taylor coefficients) at all knots.  This suggests that they will be useful for marching methods used to solve IVP for ODE \textsl{at very high precision}.  They are indeed used in this way for DAE~\cite{nedialkov2005solving}. In this section we will pursue this application.  I will end by discussing the use of the integration method for computing generalized Mathieu functions, which are needed at a double eigenvalue of the Mathieu equation.
\subsection{Collocation methods for ODE}
In this section I describe the method I implemented for~\cite{Brimacombe2021} and~\cite{Brimacombe2023}, which uses the blendstring environment described in this present paper.

One of the best codes for the numerical solution of boundary value problems (BVP) for ODE is the FORTRAN code COLSYS~\cite{ascher1979colsys}, with theoretical foundations described by the analysis in the SIAM Classic book~\cite{ascher1995numerical}.  The method that code uses is \textsl{collocation}, which requires that the differential equation be satisfied exactly at certain points, called collocation points, in each subinterval.

Their original code~\cite{ascher1979colsys} used B-splines as a piecewise polynomial interpolant, but later the authors examined Hermite and monomial bases~\cite{ascher1983spline} and found them to be superior.  Those bases were not the same as used here, because they incorporated Runge--Kutta method information as well.  There are significant simplifications that occur when, instead of Runge--Kutta information, one uses actual Taylor coefficients.

The automatic generation of Taylor series coefficients for linear differential equations with polynomial coefficients is well understood.  For instance, in Maple, one can use the \lstinline{diffeqtorec} routine from the \lstinline{gfun} package~\cite{salvy1994gfun} to generate the recurrence relations. There are other methods, some even implemented in FORTRAN~\cite{chang1994atomft}.  Generation of Taylor coefficients costs at most $O(m^2)$ floating-point operations in general, but for D-finite or holonomic systems the cost is only $O(m)$~\cite{vanderHoeven1999,van2001fast,mezzarobba2012note}.

It is efficient to use Taylor coefficients generated at \textsl{both} ends of the marching step, and blend these coefficients together.  This kind of method is called an Hermite--Obreschkov method, after the French 19th century mathematician Charles Hermite and the Bulgarian 20th century mathematician Nikola Obreschkov. Methods like this have been used with great success for differential-algebraic equations (DAE)~\cite{nedialkov2005solving} and recent work has uncovered a class of conjugate symplectic such methods~\cite{Mazzia2018}.  The main advantage to this method over an explicit Taylor series method is that by combining the terms at either end one has a method of order $2m$ instead of just $m$, at very little extra cost.

For simplicity of presentation, assume that we are solving a linear second-order scalar ODE, say $y'' +a(x)y' + b(x)y = g(x)$.  Further assume that Taylor series for the functions $a(x)$, $b(x)$, and $g(x)$ are available at working precision on demand.  What follows is one convenient and accurate method for marching from one knot to the next, by using collocation.

The method assumes that Taylor series coefficients have been generated at the current knot, say $a_n$, and are considered to be ``known".  Specifically, suppose to start with that a Taylor polynomial of grade $m$ for the desired solution is known at this knot.

Suppose also that one has chosen a tentative next knot, $a_{n+1} = a_n + h$.  If the variable $z$ were time, this would be a time step.  The stepsize $h$ is tentative at this point.  Now generate Taylor coefficients up to grade $m$ for two independent solutions, satisfying (for one solution)
\begin{equation}
    y(a_{n+1}) = 1 \ \mathrm{and}\ y'(a_{n+1}) = 0
\end{equation}
and (for the complementary solution)
\begin{equation}
    y(a_{n+1}) = 0 \ \mathrm{and}\ y'(a_{n+1}) = 1\>.
\end{equation}

Next, \textsl{blend} the known coefficients at $a_n$ with these independent solutions in the following way.
Form a blend of the known coefficients at $a_n$ with the \textsl{zero} Taylor series at $a_{n+1}$.  Call the result $L(z)$.  Form a blend of the first series above at $a_{n+1}$ with the \textsl{zero} Taylor series at $a_n$.
Call the result $C(z)$. Form a blend of the second series above with the \textsl{zero} Taylor series at $a_n$ and call the result $S(z)$.  The desired solution will then be a linear combination of these three: say $y = A\,C(z) + B\,S(z) + L(z)$.  This uses the linearity of the equation, and the linear dependence of blends on their constituent Taylor coefficients.

Now use \textsl{collocation} at the two\footnote{Two, because this example equation is second order.} points $a_n + h/4$ and $a_n + 3h/4$ (these are Chebyshev--Lobatto points, which have good properties; general such points are available for work with higher-order equations) to give us two equations in the two unknowns $A$ and $B$.  That is, compute the residuals
\begin{align}
   r_L(z) :=& \quad L'' + a(z)L' + b(z)L - g(z) \nonumber\\
   r_C(z) :=& \quad C'' + a(z)C' + b(z)C - g(z) \nonumber\\
   r_S(z) :=& \quad S'' + a(z)S' + b(z)S - g(z)
\end{align}
at those two collocation points, and set the residual for $y$ to zero at those two points:
\begin{align}
    0 = & A\> r_C(a_n + h/4) + B\> r_S(a_n+h/4) + r_L(a_n + h/4) \nonumber\\
    0 = & A\> r_C(a_n + 3h/4) + B\> r_S(a_n+3h/4) + r_L(a_n + 3h/4)\>.
\end{align}
One may solve this two-by-two linear system by any method in order to find the coefficients $A$ and $B$.  This system is nonsingular because the solutions are linearly independent at the right endpoint. Experience shows that the equations are well-scaled and well-conditioned, but this needs a proper analysis as in~\cite{bader1987new}.  One might expect that meshes with widely varying mesh sizes, or especially dense meshes, might cause problems.  However, today we have a resource that was inconvenient in the old FORTRAN codes: namely, we may use higher precision if necessary.

Now \textsl{sample} the residual (after all, we have a blend for the putative solution, and so we can evaluate it and its derivatives wherever we choose) at $a_n + h/2$, which is (asymptotically as $h \to 0$) the location of the maximum. Asymptotically as the stepsize $h \to 0$ the residual has the form
\begin{align}
r(z) =&\> y'' + a(z)y' + b(z)y - g(z) \nonumber\\
     =&\> h^{2m}Ks^{m-1}(s-1/4)(s-3/4)(s-1)^{m-1} + \mathrm{h.o.t.}\>,
\end{align}
where $s = (z-a_n)/h$, ``h.o.t" means ``higher-order terms", and $K$ is a high-order Taylor coefficient evaluated at a point between $a_n$ and $a_{n+1}$.  This truncation error is proportional to $h^{2m}$ and the maximum of the polynomial $s^{m-1}(s-1/4)(s-3/4)(s-1)^{m-1}$ occurs at $s=1/2$.  If the sampled residual is smaller than the user's tolerance, accept the step and move on; if not,  reject the step and adjust the prediction for the tentative new $a_{n+1}$ and try again.  Standard heuristics can be used here.

This presents the essence of collocation. One can see why this might be attractive in a high-precision environment.  Instead of having an $m$th order method with Taylor polynomials of grade $m$, one has a $2m$th order method. Moreover, the use of a blend gives an interpolant of high enough order to match the order of the numerical method; in contrast, a cubic spline is only suited to low-order methods.  The smallness of the error coefficients, which is determined by the infinity norm of the $s$-polynomial, makes this even more attractive.

This is an \textsl{implicit} method, which is appropriate if the ODE is \textsl{stiff}~\cite{soderlind2015stiffness}, and can help for oscillatory problems as well.
\subsection{Stability of the method for oscillatory ODE}
 When one tries this method \textsl{symbolically} on the simple harmonic oscillator $\ddot{y} + \omega^2 y = 0$, it is possible to discover an interesting stability limitation on the allowable stepsize.

To begin, exactly solve, by hand or otherwise, the simple harmonic oscillator with initial conditions $y(0) = y_0$ and $y'(0) = y_1$. The answer can be expressed as
$y(t) = y_0\cos \omega t + y_1 \sin(\omega t)/\omega $.  Now, taking a single step of size $h$ with the exact solution finds the exact solution value $Y_0$ and derivative value $Y_1$ at $t=h$, with
\begin{equation}
    \begin{bmatrix}
        Y_0 \\
        Y_1
    \end{bmatrix}
    = \begin{bmatrix}
        \cos \omega h & \frac{\sin \omega h}{\omega } \\
        -\omega \sin \omega h & \cos \omega h
    \end{bmatrix}
    \begin{bmatrix}
    y_0 \\
    y_1
    \end{bmatrix}\>.
\end{equation}
Because the equation is autonomous, this step from $t=0$ to $t=h$ is exactly the same as the same width step from $t=t_k$ to $t=t_k+h$.
Timesteps with this exact solution therefore satisfy $\mathbf{y}^{(k)} = \mathbf{A}^k \mathbf{y}^{(0)}$.
The eigenvalues of this matrix satisfy the characteristic equation $\lambda^2 -2\cos\omega h \lambda + 1 = 0$ and are $\exp(\pm i\omega h)$, which both have magnitude $1$, implying that the length of the vector of initial conditions does not grow or decay exponentially.  One would like this property to hold with the numerical method, if possible.

Applying the collocation method just described gives, at every balanced order (grade $m$ Taylor polynomials at each end), an analogous matrix, but with rational functions $C_m(\nu)$ and $S_m(\nu)$ (easily computed for any fixed $m$) of $\nu = \omega h$ in place of $\cos\nu$ and $\sin \nu$:
\begin{equation}
    \mathbf{A}_m := \begin{bmatrix}
        C_m(\nu) & \frac{S_m(\nu)}{\omega} \\
        -\omega S_m(\nu) & C_m(\nu)
    \end{bmatrix}\>.
\end{equation}
Because it turns out that $C_m^2 + S_m^2 = 1$, the matrix has characteristic polynomial $\lambda^2 - 2C_m(\nu)\lambda + 1$, implying that the product of its two eigenvalues is $1$.  The eigenvalues are $C_m(\nu) \pm iS_m(\mu)$. However, the eigenvalues \textsl{both} have magnitude $1$ if and only if $|C_m(\nu)| \le 1$.  Note that both $h$ and $\omega$, hence $\nu = \omega h$, are real.  This suggests investigating the real zeros of the equation $C_m^2(\nu)-1 = 0$. The first few rational functions $C_m(\nu)$ are tabulated in Table~\ref{tab:rats}; they are some kind of rational approximation to $\cos\nu$, but I do not recognize them (they are not $(m,m)$ Pad\'e approximants, for instance).

 \begin{table}[t]
     \centering
     \begin{tabular}{c|c|c}
     \toprule
     $m$ & $C_m(\nu)$ & $\nu^*/\pi$ \\
     \midrule
$1$          & $\frac{57 \nu^{4}-1408 \nu^{2}+3072}{9 \nu^{4}+128 \nu^{2}+3072}$ & $0.94035$\\
$2$          & $-\frac{2 \left(33 \nu^{6}-4059 \nu^{4}+84480 \nu^{2}-184320\right)}{3 \left(3 \nu^{6}+146 \nu^{4}+5120 \nu^{2}+122880\right)}$ & $0.99817$\\
$3$ & $\frac{25 \nu^{8}-9016 \nu^{6}+676560 \nu^{4}-12072960 \nu^{2}+25804800}{3 \nu^{8}+304 \nu^{6}+16080 \nu^{4}+829440 \nu^{2}+25804800}$ & $0.99997$
     \end{tabular}
     \caption{Collocation at Chebyshev--Lobatto points: The first few rational approximations to cosine and the first positive zero of $C_m^2-1$ as a fraction of $\pi$. The next entry is too wide for this table, but has $\nu^*/\pi \approx 1 - 10^{-7}$.}
     \label{tab:rats}
 \end{table}

Once the value of $C_m(\nu)$ becomes larger than $1$ in magnitude, the eigenvalues of $\mathbf{A}_m$ are no longer of unit modulus, and one of them, say $\lambda_1$, will be larger than 1 in modulus and therefore the lengths of the vectors $[y_k,y_k']$ will start to grow exponentially, like $\lambda_1^k$.  This is a numerical instability of the method.  To ensure that this does not happen, one must take $\nu < \nu^*$, or (approximately) $h < \pi/\omega$.  For high frequencies $\omega$ one would thus seem to have to take very small timesteps.

So it would seem that there is a stability restriction akin to the stepsize restrictions for explicit methods for stiff problems~\cite{soderlind2015stiffness}.  But this is not the complete story, here, and the situation is better than it seems at first: $C_m^2-1$ has \textsl{another} zero very nearby: for $m=3$, at $1.0011\pi$.  The maximum value that $C_m^2-1$ attains, on the tiny interval it is positive, is less than $3.13\cdot 10^{-6}$.  The magnitude of the largest eigenvalue is thus $1 + O(10^{-6})$. This does cause growth but, while it is technically exponential, it would not be visible in the numerical solution of the simple harmonic oscillator unless on the order of a million steps were taken!  For higher $m$, this maximum $\lambda_1$ is even smaller.  For the simple harmonic oscillator at least, this method is actually quite stable.  There are other zeros, near $2\pi$ and $3\pi$ and so on, for larger $m$, and the maxima on the small positive intervals get larger and larger until the method actually fails for large enough $\nu$, no matter how large one takes $m$.  This is because the method is not A-stable~\cite{soderlind2015stiffness}, of course.  But, A-stability is not wholly appropriate for oscillatory problems, and the current analysis gives more information.
\subsection{Computing generalized Mathieu functions}
We actually used this method, in practice~\cite{Brimacombe2021}.  We have just submitted a paper on a problem in hemodynamics where the code described in that paper is providing us with the solutions for blood flow in a tube of elliptic cross-section, using Mathieu functions.  The \textsl{reason} we did this is that no existing implementations of Mathieu functions, to our knowledge, could handle the case of double eigenvalues.  For instance, we needed to compute the Mathieu functions $\ce_0(z;q)$ and $\ce_2(z;q)$ together with the corresponding modified Mathieu functions for various purely imaginary values of the parameter $q$.  But for $q = 1.468\ldots i$, the two eigenfunctions $\ce_0(z,q)$ and $\ce_2(z,q)$ \textsl{coalesce}.  We therefore needed also to find the \textsl{generalized} eigenfunction $u(z)$ and its corresponding ``modified" generalized eigenfunction $U(z)$ in order to express the solution to our problem.  This requires solving
\begin{equation}
    u'' + (a-2q\cos 2z)u + \ce_0(z;q) = 0\>.
\end{equation}
A simple way to do this, if $\ce_0(z;q)$ is expressed as a blendstring, is to compute (on a compatible blendstring) both linearly independent solutions of the Mathieu equation: in the terminology of the DLMF (Chapter 28), these are $w_{I}(z,a,q)$ (which is really $\ce_0(z;q)$) and $w_{II}(z,a,q)$. One then forms the Green's function
\[
G(z,\zeta) = w_{I}(\zeta)w_{II}(z) - w_{II}(\zeta)w_I(z)
\]
and integrate (as blendstrings) to get
\[
u(z) = -w_{II}(z)\int_{\zeta=0}^z w_{I}(\zeta)\ce_0(\zeta)\,d\zeta  + w_I(z)\int_{\zeta=0}^z w_{II}(\zeta)\ce_0(\zeta)\,d\zeta \>.
\]
This works very well. See~\cite{Brimacombe2021}.  For the hemodynamics work using this method, which heavily uses the convenience of both the integration and differentiation of the solutions, see~\cite{Brimacombe2023}.

\begin{figure}
    \centering
    \includegraphics[width=0.305\textwidth]{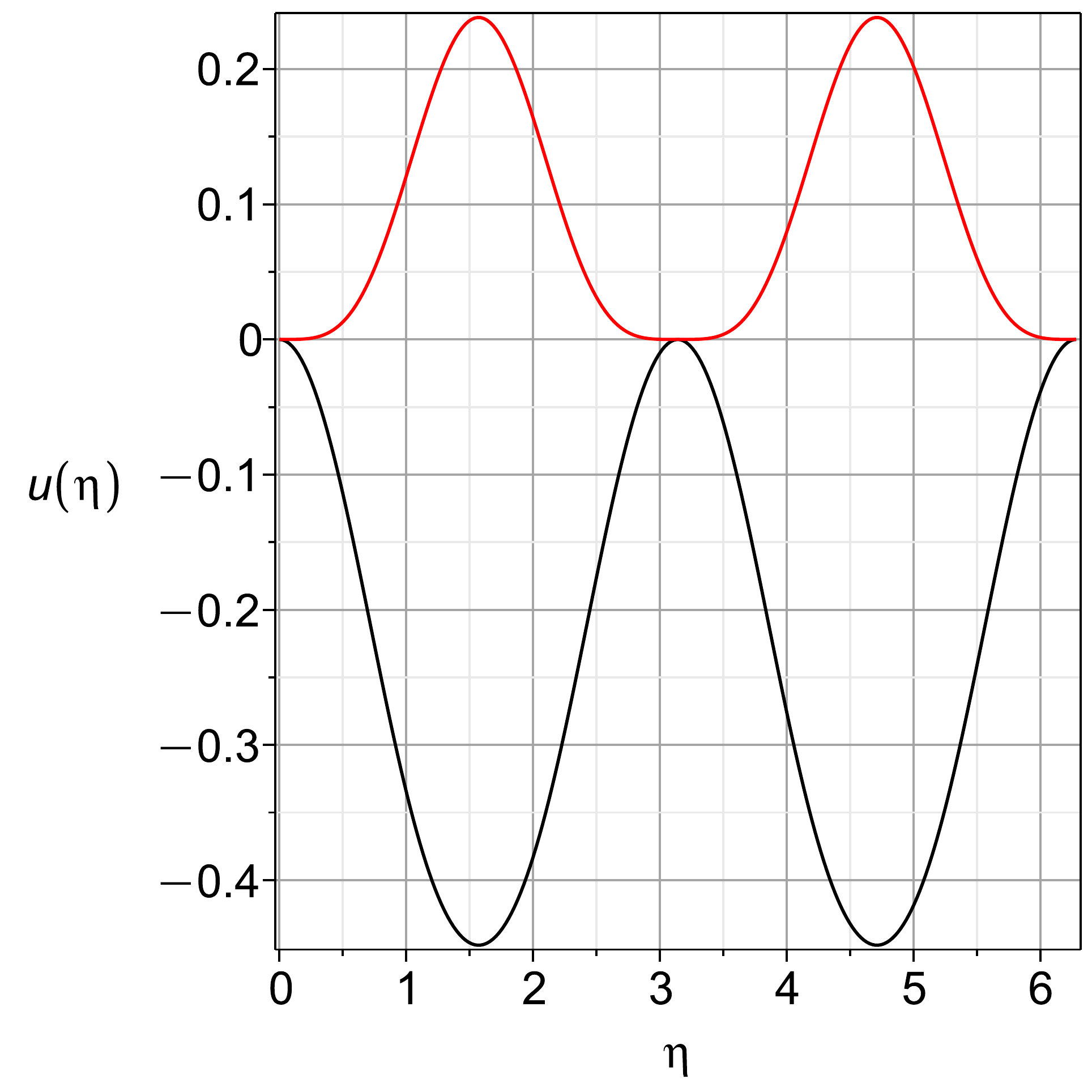}
    \caption{The generalized eigenfunction obtained when $\ce_0(x;q)$ and $\ce_2(x;q)$ coalesce at $q = 1.468\ldots i$.  Real part in black, imaginary part in red.  The solution was computed by integrating the blendstrings needed in the Green's function for the problem.}
    \Description{Two curves on a gray grid, one red (representing the imaginary part) and one black (representing the real part). Both curves are zero at the edges and the middle, which are $0$, $\pi$, and $2\pi$.}
    \label{fig:generalizedc0c2}
\end{figure}

For most ranges of the parameters, double precision suffices.  For instance when $q = 1.468\ldots\,i$, the semiminor axis $\beta = 0.151$, the semimajor axis $\alpha = 0.168$, and with grade $m=15$ Taylor polynomials (so an order $30$ method), the solver takes $395$ms to take $7$ steps across the interval\footnote{Integration could have been done on $[0,\pi]$ and saved half the time, but some of the subsequent computations concern integrals around a full period, and computation across the interval $[0,2\pi]$ made some of the bookkeeping for that simpler.} $[0,2\pi]$ to compute $\ce_0(\eta;q)$. It takes $130$ms to take $3$ steps across the (vertical) interval $[0, \xi_0\,i]$ where $\xi_0 = 1.485$ is the parameter value at the edge of the ellipse. Because of the doubly-exponential growth of $\Ce_0(\xi,q) = \ce_0(i\xi,q)$, $|\Ce_0(\xi_0,q)| = 4.7\cdot 10^8$.  Using the blendstrings and \lstinline{intBlend} to construct the generalized $u(\eta)$ shown in figure~\ref{fig:generalizedc0c2} takes no measureable further time.  All computations done in Maple 2022.1 running on a Microsoft Surface Pro 7 with Intel\textregistered\ Core\texttrademark\ i7 1065G7 1.30GHz with 4 cores and 8 logical processors.

If, however, the eccentricity $\e$ approaches $0$, the problem requires very high precision.  This seems paradoxical because circles ought to be easier than ellipses.  But the coordinate transformation used, with confocal elliptical coordinates, becomes singular here and one has $\xi \to \infty$.  Given the doubly exponential growth of Modified Mathieu functions, this becomes problematic very quickly.  Already by $\e = 0.02$ one requires hundreds of digits of precision.  We successfully used $100$ Digits and $m=80$ (giving a method of order $160$) in that case; solution still only took $4$ steps for $\ce_0$, $8$ steps for $\Ce_0$, and less than ten seconds of real time.

\section{Concluding remarks}
\begin{figure}
    \centering
    \includegraphics[width=0.305\textwidth]{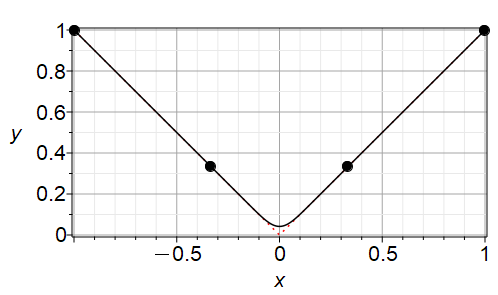}
    \caption{A blendstring on four equally-spaced knots with all Taylor polynomials of grade $m=20$, attempting to approximate $y=|x|$, which is drawn with red dashed line.  Polynomials cannot handle corners well, although at least this approximation maintains convexity.}
    \label{fig:absblend}
\end{figure}

One feature needed for the handling of derivative discontinuities is the ability to have knots with different grades of Taylor polynomials; in particular it should be possible to have just a constant (grade zero) function value at a place where the derivative does not exist. In Figure~\ref{fig:absblend} we see what happens if one tries to construct a blendstring with grade $20$ Taylor polynomials on the four equally spaced knots $-1$, $-1/3$, $1/3$ and $1$ on $[-1,1]$ for the absolute value function $f(x) = |x|$. Polynomials cannot turn sharp corners.  However, the blendstring does surprisingly well, in that it maintains convexity for this function\footnote{It works well all the way up to grade $510$ Taylor polynomials. The first visible failure (overflow, causing gaps in the graph near the knots) happens when the grade is $511$. At that point the binomials in Hermite's formula are of size $1.12\cdot 10^{306}$. Even so, the corner is visibly sharp, and perfectly placed.}.   It would be best, however, to be able to insert a knot of grade zero exactly at the location of the derivative discontinuity.  Unbalanced blends, however, have exponentially worse numerical properties, so some kind of ``multi-blend'' structure with differing series possible to the left or right of the knot might be better.

It should be possible to construct an even more accurate (and therefore efficient) kind of blend, using two-point Pad\'e approximation.
 The general global Hermite--Pad\'e construction is well-known~\cite{beckermann1992uniform}, but if the method were to be used to construct \textsl{piecewise} rational approximations with a high degree of continuity at the knots, the result may be quite interesting. Preliminary experiments are encouraging. Another generalization, to Puiseux and to Laurent series approximations on each piece, also seems worth pursuing. Use of blendstrings for matrix functions might also be interesting.
\begin{acks}
This is an outgrowth of work done jointly with Erik Postma, and with Mair Zamir \& Chris Brimacombe.  I thank Erik Postma and Michael Monagan for comments on earlier drafts.   I also thank John C.~Butcher for teaching me about the contour integral method (likely originally due to Hermite) for constructing interpolational polynomials.  After the paper was accepted, I used ChatGPT to revise my abstract to make it clearer and more concise. I then asked my friends on Mathstodon to comment on the two versions of the abstract, and took their comments into account in a synthesis of the two versions that I wrote myself.  I thank my friends for their input.
This work was supported by NSERC under RGPIN-2020-06438 and by the grant PID2020-113192GB-I00 (Mathematical Visualization: Foundations, Algorithms and Applications) from the Spanish MICINN.  I also acknowledge the support of the Rotman Institute of Philosophy.
\end{acks}




\begin{thebibliography}{23}


\ifx \showCODEN    \undefined \def \showCODEN     #1{\unskip}     \fi
\ifx \showDOI      \undefined \def \showDOI       #1{#1}\fi
\ifx \showISBNx    \undefined \def \showISBNx     #1{\unskip}     \fi
\ifx \showISBNxiii \undefined \def \showISBNxiii  #1{\unskip}     \fi
\ifx \showISSN     \undefined \def \showISSN      #1{\unskip}     \fi
\ifx \showLCCN     \undefined \def \showLCCN      #1{\unskip}     \fi
\ifx \shownote     \undefined \def \shownote      #1{#1}          \fi
\ifx \showarticletitle \undefined \def \showarticletitle #1{#1}   \fi
\ifx \showURL      \undefined \def \showURL       {\relax}        \fi
\providecommand\bibfield[2]{#2}
\providecommand\bibinfo[2]{#2}
\providecommand\natexlab[1]{#1}
\providecommand\showeprint[2][]{arXiv:#2}

\bibitem[Ascher et~al\mbox{.}(1979)]%
        {ascher1979colsys}
\bibfield{author}{\bibinfo{person}{Uri Ascher}, \bibinfo{person}{J
  Christiansen}, {and} \bibinfo{person}{Robert~D Russell}.}
  \bibinfo{year}{1979}\natexlab{}.
\newblock \showarticletitle{{COLSYS}--A collocation code for boundary-value
  problems}.
\newblock In \bibinfo{booktitle}{\emph{Codes for Boundary-Value problems in
  ordinary differential equations}}. \bibinfo{publisher}{Springer},
  \bibinfo{pages}{164--185}.
\newblock


\bibitem[Ascher et~al\mbox{.}(1995)]%
        {ascher1995numerical}
\bibfield{author}{\bibinfo{person}{Uri Ascher}, \bibinfo{person}{Robert
  Mattheij}, {and} \bibinfo{person}{Robert~D Russell}.}
  \bibinfo{year}{1995}\natexlab{}.
\newblock \bibinfo{booktitle}{\emph{Numerical solution of boundary value
  problems for ordinary differential equations}}.
\newblock \bibinfo{publisher}{SIAM}.
\newblock


\bibitem[Ascher et~al\mbox{.}(1983)]%
        {ascher1983spline}
\bibfield{author}{\bibinfo{person}{Uri Ascher}, \bibinfo{person}{Steven
  Pruess}, {and} \bibinfo{person}{Robert~D Russell}.}
  \bibinfo{year}{1983}\natexlab{}.
\newblock \showarticletitle{On spline basis selection for solving differential
  equations}.
\newblock \bibinfo{journal}{\emph{SIAM journal on Numerical Analysis}}
  \bibinfo{volume}{20}, \bibinfo{number}{1} (\bibinfo{year}{1983}),
  \bibinfo{pages}{121--142}.
\newblock


\bibitem[Bader and Ascher(1987)]%
        {bader1987new}
\bibfield{author}{\bibinfo{person}{Georg Bader} {and} \bibinfo{person}{Uri
  Ascher}.} \bibinfo{year}{1987}\natexlab{}.
\newblock \showarticletitle{A new basis implementation for a mixed order
  boundary value ODE solver}.
\newblock \bibinfo{journal}{\emph{SIAM journal on scientific and statistical
  computing}} \bibinfo{volume}{8}, \bibinfo{number}{4} (\bibinfo{year}{1987}),
  \bibinfo{pages}{483--500}.
\newblock


\bibitem[Battles and Trefethen(2004)]%
        {Battles2004}
\bibfield{author}{\bibinfo{person}{Zachary Battles} {and}
  \bibinfo{person}{Lloyd~N. Trefethen}.} \bibinfo{year}{2004}\natexlab{}.
\newblock \showarticletitle{An Extension of {MATLAB} to Continuous Functions
  and Operators}.
\newblock \bibinfo{journal}{\emph{{SIAM} Journal on Scientific Computing}}
  \bibinfo{volume}{25}, \bibinfo{number}{5} (\bibinfo{date}{Jan.}
  \bibinfo{year}{2004}), \bibinfo{pages}{1743--1770}.
\newblock
\urldef\tempurl%
\url{https://doi.org/10.1137/s1064827503430126}
\showDOI{\tempurl}


\bibitem[Beckermann and Labahn(1992)]%
        {beckermann1992uniform}
\bibfield{author}{\bibinfo{person}{Bernhard Beckermann} {and}
  \bibinfo{person}{George Labahn}.} \bibinfo{year}{1992}\natexlab{}.
\newblock \showarticletitle{A uniform approach for {H}ermite Pad{\'e} and
  simultaneous {P}ad{\'e} approximants and their matrix-type generalizations}.
\newblock \bibinfo{journal}{\emph{Numerical Algorithms}} \bibinfo{volume}{3},
  \bibinfo{number}{1} (\bibinfo{year}{1992}), \bibinfo{pages}{45--54}.
\newblock


\bibitem[Brimacombe et~al\mbox{.}(2021)]%
        {Brimacombe2021}
\bibfield{author}{\bibinfo{person}{Chris Brimacombe},
  \bibinfo{person}{Robert~M. Corless}, {and} \bibinfo{person}{Mair Zamir}.}
  \bibinfo{year}{2021}\natexlab{}.
\newblock \showarticletitle{Computation and applications of {M}athieu
  functions: A historical perspective}.
\newblock \bibinfo{journal}{\emph{SIAM Rev.}} \bibinfo{volume}{63},
  \bibinfo{number}{4} (\bibinfo{date}{Jan.} \bibinfo{year}{2021}),
  \bibinfo{pages}{653--720}.
\newblock
\urldef\tempurl%
\url{https://doi.org/10.1137/20m135786x}
\showDOI{\tempurl}


\bibitem[Brimacombe et~al\mbox{.}(2023)]%
        {Brimacombe2023}
\bibfield{author}{\bibinfo{person}{Chris Brimacombe},
  \bibinfo{person}{Robert~M. Corless}, {and} \bibinfo{person}{Mair Zamir}.}
  \bibinfo{year}{2023}\natexlab{}.
\newblock \showarticletitle{Elliptic cross sections in blood flow regulation}.
\newblock \bibinfo{journal}{\emph{ArXiv}} (\bibinfo{year}{2023}).
\newblock
\urldef\tempurl%
\url{https://arxiv.org/abs/2304.01356}
\showURL{%
\tempurl}


\bibitem[Chang and Corliss(1994)]%
        {chang1994atomft}
\bibfield{author}{\bibinfo{person}{YF Chang} {and} \bibinfo{person}{George
  Corliss}.} \bibinfo{year}{1994}\natexlab{}.
\newblock \showarticletitle{{ATOMFT}: solving {ODE}s and {DAE}s using {T}aylor
  series}.
\newblock \bibinfo{journal}{\emph{Computers \& Mathematics with Applications}}
  \bibinfo{volume}{28}, \bibinfo{number}{10-12} (\bibinfo{year}{1994}),
  \bibinfo{pages}{209--233}.
\newblock


\bibitem[Corless(2023)]%
        {Corless2023details}
\bibfield{author}{\bibinfo{person}{Robert~M. Corless}.}
  \bibinfo{year}{2023}\natexlab{}.
\newblock \showarticletitle{Blends have decent numerical properties}.
\newblock \bibinfo{journal}{\emph{Maple Transactions}} \bibinfo{volume}{3},
  \bibinfo{number}{1} (\bibinfo{date}{Feb.} \bibinfo{year}{2023}).
\newblock
\urldef\tempurl%
\url{https://doi.org/10.5206/mt.v3i1.15890}
\showDOI{\tempurl}


\bibitem[Corless and Postma(2021)]%
        {Corless2021}
\bibfield{author}{\bibinfo{person}{Robert~M. Corless} {and}
  \bibinfo{person}{Erik~J. Postma}.} \bibinfo{year}{2021}\natexlab{}.
\newblock \showarticletitle{Blends in {Maple}}.
\newblock In \bibinfo{booktitle}{\emph{Communications in Computer and
  Information Science}}. \bibinfo{publisher}{Springer International
  Publishing}, \bibinfo{pages}{167--184}.
\newblock
\urldef\tempurl%
\url{https://doi.org/10.1007/978-3-030-81698-8_12}
\showDOI{\tempurl}


\bibitem[Geddes and Fee(1992)]%
        {geddes1992hybrid}
\bibfield{author}{\bibinfo{person}{Keith~O Geddes} {and}
  \bibinfo{person}{Gregory~J Fee}.} \bibinfo{year}{1992}\natexlab{}.
\newblock \showarticletitle{Hybrid symbolic-numeric integration in {M}APLE}. In
  \bibinfo{booktitle}{\emph{Papers from the international symposium on Symbolic
  and algebraic computation}}. \bibinfo{pages}{36--41}.
\newblock


\bibitem[Hermite(1873)]%
        {hermite1873cours}
\bibfield{author}{\bibinfo{person}{Charles Hermite}.}
  \bibinfo{year}{1873}\natexlab{}.
\newblock \bibinfo{booktitle}{\emph{Cours d'analyse de l'{\'E}cole
  polytechnique}}. Vol.~\bibinfo{volume}{25}.
\newblock \bibinfo{publisher}{Gauthier-Villars}.
\newblock


\bibitem[Ilie et~al\mbox{.}(2008)]%
        {ilie2008adaptivity}
\bibfield{author}{\bibinfo{person}{Silvana Ilie}, \bibinfo{person}{Gustaf
  S{\"o}derlind}, {and} \bibinfo{person}{Robert~M Corless}.}
  \bibinfo{year}{2008}\natexlab{}.
\newblock \showarticletitle{Adaptivity and computational complexity in the
  numerical solution of ODEs}.
\newblock \bibinfo{journal}{\emph{Journal of Complexity}} \bibinfo{volume}{24},
  \bibinfo{number}{3} (\bibinfo{year}{2008}), \bibinfo{pages}{341--361}.
\newblock


\bibitem[Mazzia and Sestini(2018)]%
        {Mazzia2018}
\bibfield{author}{\bibinfo{person}{Francesca Mazzia} {and}
  \bibinfo{person}{Alessandra Sestini}.} \bibinfo{year}{2018}\natexlab{}.
\newblock \showarticletitle{On a Class of Conjugate Symplectic
  {Hermite--Obreshkov} One-Step Methods with Continuous Spline Extension}.
\newblock \bibinfo{journal}{\emph{Axioms}} \bibinfo{volume}{7},
  \bibinfo{number}{3} (\bibinfo{date}{Aug.} \bibinfo{year}{2018}),
  \bibinfo{pages}{58}.
\newblock
\urldef\tempurl%
\url{https://doi.org/10.3390/axioms7030058}
\showDOI{\tempurl}


\bibitem[Mezzarobba(2012)]%
        {mezzarobba2012note}
\bibfield{author}{\bibinfo{person}{Marc Mezzarobba}.}
  \bibinfo{year}{2012}\natexlab{}.
\newblock \showarticletitle{A note on the space complexity of fast {D}-finite
  function evaluation}. In \bibinfo{booktitle}{\emph{Int. Workshop on Computer
  Algebra in Scientific Computing}}. Springer, \bibinfo{pages}{212--223}.
\newblock


\bibitem[Nedialkov and Pryce(2005)]%
        {nedialkov2005solving}
\bibfield{author}{\bibinfo{person}{Nedialko~S Nedialkov} {and}
  \bibinfo{person}{John~D Pryce}.} \bibinfo{year}{2005}\natexlab{}.
\newblock \showarticletitle{Solving differential-algebraic equations by
  {Taylor} series (I): Computing {Taylor} coefficients}.
\newblock \bibinfo{journal}{\emph{BIT Numerical Mathematics}}
  \bibinfo{volume}{45}, \bibinfo{number}{3} (\bibinfo{year}{2005}),
  \bibinfo{pages}{561--591}.
\newblock


\bibitem[Salvy and Zimmermann(1994)]%
        {salvy1994gfun}
\bibfield{author}{\bibinfo{person}{Bruno Salvy} {and} \bibinfo{person}{Paul
  Zimmermann}.} \bibinfo{year}{1994}\natexlab{}.
\newblock \showarticletitle{Gfun: a {M}aple package for the manipulation of
  generating and holonomic functions in one variable}.
\newblock \bibinfo{journal}{\emph{ACM Transactions on Mathematical Software
  (TOMS)}} \bibinfo{volume}{20}, \bibinfo{number}{2} (\bibinfo{year}{1994}),
  \bibinfo{pages}{163--177}.
\newblock


\bibitem[Shampine and Reichelt(1997)]%
        {shampine1997matlab}
\bibfield{author}{\bibinfo{person}{Lawrence~F Shampine} {and}
  \bibinfo{person}{Mark~W Reichelt}.} \bibinfo{year}{1997}\natexlab{}.
\newblock \showarticletitle{The {MATLAB} {ODE} suite}.
\newblock \bibinfo{journal}{\emph{SIAM journal on scientific computing}}
  \bibinfo{volume}{18}, \bibinfo{number}{1} (\bibinfo{year}{1997}),
  \bibinfo{pages}{1--22}.
\newblock


\bibitem[S{\"o}derlind et~al\mbox{.}(2015)]%
        {soderlind2015stiffness}
\bibfield{author}{\bibinfo{person}{Gustaf S{\"o}derlind},
  \bibinfo{person}{Laurent Jay}, {and} \bibinfo{person}{Manuel Calvo}.}
  \bibinfo{year}{2015}\natexlab{}.
\newblock \showarticletitle{Stiffness 1952--2012: Sixty years in search of a
  definition}.
\newblock \bibinfo{journal}{\emph{BIT Numerical Mathematics}}
  \bibinfo{volume}{55}, \bibinfo{number}{2} (\bibinfo{year}{2015}),
  \bibinfo{pages}{531--558}.
\newblock


\bibitem[Trefethen(2019)]%
        {trefethen2019approximation}
\bibfield{author}{\bibinfo{person}{Lloyd~N Trefethen}.}
  \bibinfo{year}{2019}\natexlab{}.
\newblock \bibinfo{booktitle}{\emph{Approximation Theory and Approximation
  Practice}}.
\newblock \bibinfo{publisher}{SIAM}.
\newblock


\bibitem[van~der Hoeven(1999)]%
        {vanderHoeven1999}
\bibfield{author}{\bibinfo{person}{Joris van~der Hoeven}.}
  \bibinfo{year}{1999}\natexlab{}.
\newblock \showarticletitle{Fast evaluation of holonomic functions}.
\newblock \bibinfo{journal}{\emph{Theoretical Computer Science}}
  \bibinfo{volume}{210}, \bibinfo{number}{1} (\bibinfo{date}{Jan.}
  \bibinfo{year}{1999}), \bibinfo{pages}{199--215}.
\newblock
\urldef\tempurl%
\url{https://doi.org/10.1016/s0304-3975(98)00102-9}
\showDOI{\tempurl}


\bibitem[van~der Hoeven(2001)]%
        {van2001fast}
\bibfield{author}{\bibinfo{person}{Joris van~der Hoeven}.}
  \bibinfo{year}{2001}\natexlab{}.
\newblock \showarticletitle{Fast evaluation of holonomic functions near and in
  regular singularities}.
\newblock \bibinfo{journal}{\emph{Journal of Symbolic Computation}}
  \bibinfo{volume}{31}, \bibinfo{number}{6} (\bibinfo{year}{2001}),
  \bibinfo{pages}{717--744}.
\newblock


\end{thebibliography}
\end{document}